\documentclass[12pt]{article}
\usepackage{times}
\usepackage{a4wide}
\usepackage{amsfonts}
\usepackage{theorem}
\usepackage{enumerate}
\usepackage{amsmath}
\usepackage{amscd}
\usepackage{amssymb}

\theorembodyfont{\sl}

\newtheorem{thm}[subsection]{Theorem}
\newtheorem{prop}[subsection]{Proposition}
\newtheorem{lem}[subsection]{Lemma}
\newtheorem{cor}[subsection]{Corollary}
\newtheorem{sublem}[subsubsection]{Sublemma}

\theorembodyfont{\rmfamily}

\newtheorem{rem}[subsection]{Remark}









\newenvironment{prf}[1]{\trivlist
\item[\hskip \labelsep{\it
#1.\hspace*{.3em}}]}{~\hspace{\fill}~$\square$\endtrivlist}
\newenvironment{proof}{\begin{prf}{\bf Proof}}{\end{prf}}

\newcommand{\CC}{{\mathbb C}}
\newcommand{\RR}{{\mathbb R}}

\newcommand{\QQ}{{\mathbb Q}}
\newcommand{\FF}{{\mathbb F}}
\newcommand{\ZZ}{{\mathbb Z}}
\newcommand{\PP}{{\mathbb P}}

\newcommand{\GG}{{\mathbb G}}
\newcommand{\TT}{{\mathbb T}}

\newcommand{\EE}{{\mathbb E}}

\newcommand{\Spec}{{\rm Spec}}

\newcommand{\Pic}{{\rm Pic}}

\newcommand{\into}{\hookrightarrow}

\newcommand{\calV}{{\cal V}}
\newcommand{\calO}{{\cal O}}
\newcommand{\calE}{{\cal E}}
\newcommand{\calC}{{\cal C}}

\newcommand{\calG}{{\cal G}}

\newcommand{\bs}{\backslash}

\newcommand{\ol}{\overline}

\newcommand{\End}{{\rm End}}

\newcommand{\GL}{{\rm GL}}
\newcommand{\Frob}{{\rm Frob}}

\title{Proving the triviality of rational points on Atkin-Lehner quotients of 
Shimura curves}
\date{}
\author{Pierre Parent and Andrei Yafaev}

\begin{document}
\maketitle

\begin{abstract}
In this paper we give a method for studying global rational points on 
certain quotients of Shimura curves by Atkin-Lehner involutions. We 
obtain explicit conditions on such quotients for rational points to be
``trivial'' (coming from CM points only) and exhibit an explicit  
infinite family of such quotients satisfying these conditions.
\end{abstract}

\section{Introduction.}

The study of points on Shimura curves rational over number fields or,
more precisely, finding families of such curves whose only rational
points are special (or CM) is an interesting problem for several reasons. One 
reason is that these curves may provide non-trivial examples of violation of 
the Hasse principle accounted for by the Brauer-Manin obstruction. This 
question was the main motivation for \cite{SkoYaf}, \cite{RotSkoYaf}, 
\cite{Sko} and \cite{SkoSiksek}.

Another motivation is the question of scarcity of possible endomorphism 
algebras of abelian varieties over $\QQ$ of $\GL_2$-type. Let us explain this in 
more detail, following the introduction to \cite{BruiFlyGonRot}. One says 
that an abelian variety $A$ over $\QQ$ is of $\GL_2$-type if the endomorphism 
algebra $\QQ\otimes \End_{\QQ}(A)$ is a number field of degree $\dim(A)$. 
In loc. cit. the authors formulate the following conjecture (attributed to 
Coleman, also made by Clark and Mazur): for $g$ any positive integer, 
the number of isomorphism classes of $\overline{\QQ}$-endomorphism 
algebras of $g$-dimensional abelian varieties of $\GL_2$-type is 
finite. We refer to loc. cit. for more details on this and just mention that 
the ``$\GL_{2}$-type'' hypothesis in the conjecture is motivated by a 
conjecture of Ribet predicting that abelian varieties over $\QQ$ of 
$\GL_2$-type are modular: in the latter case, the extended experimental 
knowledge available on modular abelian varieties gives evidence supporting 
the conjecture. Note also that, as Ribet himself has shown, his  
conjecture follows from Serre's conjecture on modularity of mod $p$
representations, and that a proof of this conjecture has recently been 
announced by Khare et al.  For $g=1$, the theory of complex 
multiplication, plus the classical result that there is only a finite number 
of quadratic imaginary orders with class number 1, show that Coleman's
conjecture is true. For $g\geq 2$ very little is known. When $g=2$ the 
endomorphism algebras that may occur are orders in either $M_2 (\QQ)$, 
or $M_2 (K)$ for $K$ an imaginary quadratic field (in these cases
the abelian surfaces are geometrically reducible), or (in the irreducible 
case) a quadratic field or an indefinite division quaternion algebra over 
$\QQ$ (see \cite{BruiFlyGonRot}, introduction and references therein). 

The diophantine study of Shimura curves allows to address the
question of excluding this latter case of (potential) quaternionic
multiplication, as these curves are coarse moduli spaces for abelian
surfaces endowed with a geometric action of some order in a quaternion
algebra. Indeed if $D$ is an integer which is the product of an even
number of distinct primes and $m$ is a square-free integer, we say that
the pair $(D,m)$ is of $\GL_2$-type if there is an abelian surface $A$
over $\QQ$ such that $\End_{\QQ} (A)\otimes \QQ =\QQ (\sqrt{m} )$ 
and $\End_{\overline{\QQ}} (A)\otimes \QQ =B_D$ (the quaternion 
algebra ramified at primes dividing $D$). Let $X^D$ be the Shimura 
curve over $\QQ$ attached to a maximal order in $B_D$ (denoted by 
$X_D$ in \cite{BruiFlyGonRot}), and if $n$ is a divisor of $D$, let 
$w_n$ be the Atkin-Lehner involution associated to $n$. Recall that a 
point on these curves is said to be a special (or CM) point if the 
underlying abelian surface $A$ is isogenous to a square of an elliptic 
curve having complex multiplication by a field $K$ (in which case  
${\mathrm{End}}_{\overline{\QQ}} (A)\otimes \QQ \simeq M_2
(K)$ contains any rational quaternion algebra split over $K$). With the 
above notations Rotger proves the following result: if $(D,m)$ is a pair 
of $\GL_2$-type over $\QQ$, then $m|D$ and $X^D /w_m (\QQ )$ is 
{\em non-trivial}, i.e. contains a non-special point (see 
\cite{RotgerTransAMS}, Theorem 4.4).

It has long been known that $X^D(\QQ)$ is empty (in fact $X^D(\RR)
=\emptyset$). More generally the results of Jordan, Livn\'e, Ogg, 
Shimura and some others give a description of the set of {\it local} 
points on Shimura curves and their Atkin-Lehner quotients. Their 
results provide necessary conditions for  $X^D /w_m$ to have local
points everywhere, and this leads to the following (see \cite{RotSkoYaf}, 
Theorem 3.1):
\begin{thm}
If $X^D /w_m (\QQ ) \neq \emptyset$, then either $m=D$ or
$m=D/p$ for some $p$ dividing $D$. Furthermore if $(D,m)=(pq,q)$ for
two odd primes $p$ and $q$ then $\left( \frac{q}{p} \right) =-1$ and
either:
\begin{enumerate}
\item $p\equiv 3\mod 4$, or
\item $p\equiv 1 \mod 4$ and $q\equiv 3 \mod 4$.
\end{enumerate}
\end{thm}

The case $m=D$ is one of the main questions addressed by P. Clark in
\cite{Clark}. Unfortunately, the method we develop in this paper does 
not work in this case. The reason is that our method relies, in a crucial 
way, on the existence of a non-trivial rank zero quotient of the jacobian
of the curve. Considerations of the sign of the functional equation of 
$L$-functions show that the jacobian of the curve $X^D/w_D$ has no 
non-trivial rank-zero quotient (at least if one believes in the Birch and 
Swinnerton-Dyer conjecture).

However in the case where $m=D/p$, there is no obstruction to the
existence of such a quotient coming from the sign of the functional equation. 
Let $X$ be the curve $X^D/w_m$. We are interested in determining whether 
$X(\QQ)$ is ``trivial'' (in the sense that it consists only of special points). Our 
strategy follows the lines of the paper \cite{ParentX0+} of the first author, 
with several new ingredients. More precisely a local study as in 
\cite{OggMR}, section 1, leads to considering two subcases that Ogg  
calls ``ramifi\'e'' and ``non-ramifi\'e'' (see loc. cit., p. 204 et seq.) that
correspond respectively to the two cases at the end of the above theorem. 
For technical reasons it is easier to push our methods further in the latter 
(``non-ramifi\'e'') case, and we limit ourselves to this case in this paper. 
For a given pair $(pq,q)$ one can in principle verify, by computation, whether 
our conditions are satisfied. However in this article we were not able to make 
our method work for infinitely many $q$s. We nevertheless conclude our paper 
by exhibiting an infinite family of Atkin-Lehner quotients of Shimura curves 
whose only possible rational points are trivial in the above sense, although 
these curves may have local points at every prime of bad reduction. (Actually, 
one can also check that in our family, there are no $\QQ$-rational special 
points either, so the set of rational points is empty; see section 6). As far as 
we know, this was previously known for finitely many curves only (but see 
also Rotger's result discussed below).

   We close this introduction by discussing a recent result, similar 
to the one obtained in this article. After the present paper was completed we 
learned from V. Rotger that he could also prove the triviality of rational points
in infinite families of Atkin-Lehner quotients of Shimura curves, but in the 
``ramifi\'e" case quoted above, more precisely for pairs $(pq,q)$ as in the first 
case of Theorem 1.1 (see \cite{Rotger}, Theorem 1.2). His methods are 
totally different from ours: he uses the modular interpretation of the points 
on the Shimura curves and analyzes the Galois representations provided by 
torsion points of the corresponding abelian surfaces. Note that our results 
are complementary as we focus on the  ``non-ramifi\'e" case. In fact, as 
Rotger himself pointed out to us, it follows from his work that in this latter 
case the field of definition of modular points is different from the field of 
moduli, i.e. the $\QQ$-valued points correspond to abelian surfaces 
that {\it cannot} be defined over $\QQ$. On the other hand, our 
methods do not use the modular interpretation of the points, so they
cover cases that cannot be tackled by Rotger's techniques.
\subsection{The strategy.}
We now sketch our method. Suppose that $D$ is a product of two primes
$D=pq$ and $m$ is $q$. Let $X:=X^{pq} /w_q$ be as above and suppose that 
there is a point $P$ in $X(\QQ )$. Assume that we are in the case 2 of
Theorem 1.1. Consider the model of $X$ over $\ZZ_p$ given
by the Cherednik-Drinfeld theorem and let $\tilde{X}$ be the regular 
model of $X$ obtained by blowing-up at singular points. Looking at the
specialization of $P$ in $\tilde{X}$ at $p$, one sees as in \cite{OggMR},
p. 205, that $P$ specializes precisely to an exceptional component
coming from a blow-up of a singular point of width 2 in the special
fiber of $X$, i.e. a singular point having ``multiplication by $i$''. The same 
is true for $w_p (P)$, therefore the specializations of our two rational
points $P$ and $w_p(P)$ belong to the smooth locus of the {\it same}
exceptional component at $p$. Let $\phi_{P}$ be the Albanese embedding 
of $X^{\mathrm{sm}}_{/\ZZ_p}$ into the N\'eron model of its jacobian 
$J$ defined by $\phi_{P} (Q)=(Q)-(P)$. Let $X_{p+1,\QQ}$ be the 
curve $(p+1)\phi_{P} (X_\QQ)$ in $J_\QQ$.

Assuming that $p$ is much larger than $q$ we prove that $X_{p+1,\QQ}$ 
has, at least locally, a smooth model around $0_{\FF_p}$, and the points 
$(p+1)\phi_{P} (P)$ and $(p+1)\phi_{P} (w_p (P))$ coincide in $X_{p+1} 
(\FF_p )$. This is the content of sections 2 and 3. We hope that the 
description we give here of the asymptotic shape (as the discriminant 
increases) of the graph of special fibers of Shimura curves, and therefore of 
the component groups of their jacobians, can be of independent interest. 

Next we consider the natural map $\Psi_{P}$ from $X$ to the ``winding 
quotient'' $J_e$ of $J$. We prove as in \cite{ParentX0+}
that, if one can find on the graph of $X$ a closed path ``made'' of
Gross-Heegner vectors containing the component of $P(\FF_p )$, then 
$\Psi_{P}$ is a formal immersion at our point $P_{\FF_p}$. This is proved 
in section 5, using the results of section 4. As $p+1$ is prime to $p$, the
same is true for the map from $X_{p+1}$ to $J_e$ at $(p+1)\phi_{P}
(P)_{\FF_p}$. The classical argument of Mazur's method now applies, 
giving that $(p+1)\phi_{P} (P)=(p+1)\phi_{P} (w_{p} (P))$ in $X_{p+1}
(\QQ )$, so $(p+1)((P)-(w_p (P)))=0$ in $J(\QQ )$. This means that either 
$P=w_p (P)$, in which case $P$ is indeed a special point, or $X$ has 
gonality less than $p+1$. If $q$ is larger than 245 this contradicts 
Abramovich's bound on the gonality of Shimura curves, and this argument 
concludes the proof (see Theorem 5.3).

One sees from the above sketch that one of the main problems
is to construct a ``closed path made of Gross-Heegner vectors
containing the component of $P(\FF_p )$ on the graph of $X$''. This is what
we do in the last section for a fixed $q=251$, and $p$ ranging through an
infinite set of primes satisfying certain congruences.

{\it Acknowledgments.}\\
We are grateful to Bas Edixhoven, Qing Liu, Victor Rotger and Emmanuel
Ullmo for instructive discussions on the subject of this paper.
We thank the referee for useful comments and suggestions.
\\
\section{Admissible curves and their jacobians.}
We start by recalling the notion of {\it admissible curve} in the
sense of Jordan and Livn\'e, and some of their basic properties.
The theorem of Cherednik and Drinfeld implies that
the Shimura curves and their Atkin-Lehner quotients
have natural models which are admissible curves at primes dividing
their discriminant. For more details we refer to \cite{JordanLivne},
section 3, or \cite{Ribet}, section 2.

Let $p$ be a prime number and $K$ a $p$-adic field with ring of
integers ${\cal O}$, uniformizer $\pi$ and residue field $k$. Let $X$
be a proper and flat curve over ${\cal O}$. We say that $X$ is
admissible if its generic fiber $X_K$ is smooth over $K$, its special
fiber $X_k$ is potentially (i.e., after a finite extension of $k$) a
union of projective lines with reduced normal crossings, and each
singularity of $X_k$ has a local equation $xy=\pi^e$ where $e\geq 1$ is
some integer called the width of the singularity. The curve $X$ is
regular at a singular point if and only if the width of this point is
$1$. From an arbitrary admissible curve $X$ one obtains a regular one
$\tilde{X}$ by making successive blow-ups at the singular points.
This leads to replacing a singularity of width $e$ by a chain of $e-1$
projective lines whose intersections points have width $1$. We write 
$X^{\mathrm{sm}}$ for the largest open subset of $X$ which is smooth
over ${\cal O}$.

The fiber $X_k$ of an admissible curve is combinatorially
described by its dual graph, which we denote by ${\cal G}(X_k)$.
Recall that this is the unoriented graph whose vertices
$\cal V$ are the components of $X_k$, and whose edges $\cal E$
are the singular points of $X_k$: an edge links two vertices if the
corresponding singular point is the intersection point of the two
corresponding components. One endows each edge with a length, 
which is defined to be the width of the corresponding singularity.

The neutral component $\Pic^0 (X_k )$ of the Picard scheme of $X_k$
is a torus whose character group is isomorphic to the homology group 
$H_1 ({\cal G}(X_k ), \ZZ )$. It follows from the above that the
regularization of $X_k$ amounts to replacing each edge of length $e$ in 
${\cal G} (X_k )$ by a chain of $e$ edges of length $1$. This shows 
that $H_1 ({\cal G}(X_k ), \ZZ )$ and $H_1 ({\cal G} 
(\tilde{X}_k ), \ZZ )$ are canonically isomorphic. There is a canonical 
Euclidean pairing (the monodromy pairing) on $H_1 ({\cal G}(X_k ), 
\ZZ )$ which is deduced from the pairing on $\ZZ^{\cal E}$ defined by 
$\langle e_i ,e_j \rangle =\mathrm{width} (e_i )^{-1} \delta_{i,j}$ 
(here $\delta$ denotes the Kronecker symbol). The following 
proposition explicits the isomorphism $\Pic^0 (X_k )\simeq H_1 
({\cal G} (X_k ), \ZZ ) \otimes \GG_m$ of \cite{SGA7}, 12.3.7.
\begin{prop}
Let $X$ be an admissible curve over $\calO$ as above and let
$X_{\overline{k}}$ be its special fibre. Let $C_0$ be a component of
$X_{\overline{k}}$, and $H$ be the subgroup of
$\Pic^0(X_{\overline{k}} )$ consisting of linear equivalence classes
of Weil divisors of degree zero on
$X_{\overline{k}}$ whose support is contained in the smooth locus
$C_0^{\mathrm{sm}}$ of $C_0$. Let $\{{\calC}_i \}$, $i=1,\dots, n$ be
the set of connected components of $X_{\overline{k}} \bs
C_0^{\mathrm{sm}}$. For each $i$ let $A_{ij}$, $j=1,\dots, n_i$ be the
set of singular points of $C_0 \cap \calC_i$. Let $x$ be a coordinate
function on the rational curve $C_0$.

There is a $\overline{k}$-isomorphism from $\Pic^0 (X_{\overline{k}} )$
to $T_1 \times \dots \times T_n \times {\cal T}$, where the $T_i$'s
and $\cal T$ are tori with $\mathrm{dim} (T_i )=n_i -1$, such that any
element $\sum_r (P_r )-\sum_r (Q_r )$ of $H$ is mapped to $((\prod_r
\frac{x(A_{i1})-x(Q_r )}{x(A_{ij})-x(Q_r )} \cdot \frac{x
(A_{ij})-x(P_r )}{x(A_{i1})-x(P_r )})_{2\leq j \leq n_i} )_{1\leq i
\leq n} \times 1$ in $(\prod_{i=1}^n T_i (\overline{k} ))\times {\cal
T} (\overline{k} )$.
\end{prop}
\begin{proof}
Let $D$ be a divisor on $X_{\overline{k}}$. Let $f_D$ be a function on
the normalization of $X_{\overline{k}}$ having divisor $D_{|C}$ on
each component $C(\simeq \PP^1_{\overline{k}} )$. Let $\cal B$ be a
basis of $H_1 ({\cal G} (X_{\overline{k}} ),\ZZ )$. One associates to 
$D$ the element $(\prod_{h_i , h=\sum h_i } (f_D (s(h_i ))/f_D (t(h_i)) )
)_{h\in {\cal B}}$ of $H_1 ({\cal G} (X_{\overline{k}} ), \ZZ )
\otimes \GG_m$, where the symbol $h=\sum h_i$ denotes the path $h$ 
consisting of successive edges $h_i$, and $s(h_i )$, $t(h_i )$ denote the points 
of $X_{\overline{k}}$ corresponding to the source and target of $h_i$ 
respectively (with respect to some fixed orientation on the graph).

As $H_1 ({\cal G}(X_{\overline{k}} ), \ZZ )=\oplus_{1 \leq i\leq n}
H_1 ({\cal G} ({\calC}_i \cup C_0 ), \ZZ )$, we can suppose without loss
of generality that $n=1$. We have $H_1(\calG (X_{\overline{k}}
),\ZZ)= H_1(\calG',\ZZ) \oplus H_1 (\calG ({\calC}_1 ),\ZZ) $ where
$\calG'$ is the graph of some sub-curve $C'$ of $X_{\overline{k}}$,
and this splitting induces a decomposition $\Pic^0 ({X_{\overline{k}} } )
\simeq H_1 (\calG',\ZZ) \otimes \GG_m \times H_1(\calG ({\calC}_1
),\ZZ) \otimes \GG_m$. The subgroup $H$ is mapped isomorphically to the
first factor, and has trivial image in the second. Indeed, as above, we
associate to an element $D=\sum_{r} (P_r )-(Q_r )$ of $H$ a function $f_D$
on $X_{\overline{k}}$ defined by $f_D \equiv 1$ on the complement of
$C_0$ and $f_D (x)=\prod_r \frac{x-x(P_r )}{x-x(Q_r )}$ on $C_0$. The 
map sending $D$ to $(f_D (A_{1,j}) /f_D (A_{1,1} ))_{2\leq j\leq n_1}$
identifies $H$ with $\GG_m^{n_1 -1}$.
\end{proof}

This proposition describes the Albanese images of an admissible curve in 
its jacobian. We are now interested in models of the {\it multiples} of 
these images. More precisely, letting $X$ be an admissible curve over 
$\Spec ({\cal O})$, we let $J_K :={\mathrm{Jac}} (X_K )$ and $J$ be 
its N\'eron model over $\Spec ({\cal O})$. Let $P$ and $P_0$ be two 
points of $X({\cal O})$. Call $\phi \colon 
X^{\mathrm{sm}} \to J$ the Albanese embedding defined by $\phi (Q)
:=(Q)-(P_0 )$. For any integer $n$ the image of $X_K$ by $n\cdot \phi_K$ 
in the generic fiber of $J$ is a proper curve over $K$ which we call 
$X_{n,K}$. Following \cite{BasRamanujan}, section 9, we say that a point 
on an admissible curve is {\it non-disconnecting} if the special fiber of 
the curve deprived from this point remains connected.
\begin{prop}
Assume the following:
\begin{enumerate}
\item the points $P$ and $P_0$ specialize to the same component $C_0$
of $X_k$;
\item the component $C_0$ intersects other components at two singular
points only, which are non-disconnecting;
\item $\phi (X_K )\cap J_K [n]=\{ 0\}$;
\item the image of $X$ by $\phi$ has no non-trivial $n$-torsion in the
component group of $J$, i.e. if $C\neq C_0$ is a component of $X_k$
then $n\cdot \phi (C^{\mathrm{sm}} )\not\subseteq J^0_k$.
\end{enumerate}
Then there is a smooth model of $X_{n,K}$ over $\cal O$, i.e. a smooth
$\Spec ({\cal O})$-scheme $X_n$ of relative dimension 1 such that
$X_n \times_{\Spec ({\cal O})} \Spec (K)=X_{n,K}$, which contains 
$n\cdot \phi (P)(\Spec ({\cal O}))$ and $n\cdot \phi (P_{0})(\Spec 
({\cal O}))$.
\end{prop}
\begin{proof} Let $J^0$ be the neutral component of $J$ and let
$\Pic^0 (X)$ be the usual functor from the category of $\Spec ({\cal
O})$-schemes to the category of groups which associates to $S$ the
group of invertible sheaves having degree $0$ on every component of
$X\times_{\Spec ({\cal O})} S$. As $X$ is semistable over ${\cal O}$
the scheme $J^0$ is semiabelian, so it represents $\Pic^0 (X)$ (see
\cite{RayJac}, proof of Th\'eor\`eme 2, $i$)). Therefore the assumption 
1 implies that $J^0$ contains $\phi (P)$ and $\phi (P_0 )=0$. Let $X_n$ 
be the scheme-theoretic closure of $n\cdot \phi (X)_K$ in $J^0$, i.e. the 
smallest closed subscheme (with reduced induced structure) of $J^0$ 
containing $n\cdot \phi (X)_K$. Let $A_1$ and $A_2$ be the singular 
points of $C_0$. Proposition 2.1 says that for $Q\in C_0 (\overline{k}) 
\backslash \{ A_1 , A_2 \}$,
$$\phi (Q)=\lambda (x(Q)-x(A_1))/(x(Q)-x(A_2)) \times 1$$
in a suitable decomposition $\Pic^0 (X_{\overline{k}} )\simeq T_1 
\times {\cal T}$ (with $\lambda =(x(P_0 )-x(A_2))/(x(P_0 )-x(A_1))$). 
One deduces from that and the assumption 2 that the special fibre of $X_n$ 
contains the same one-dimensional subtorus $T_1$ of $J^0_k$ as the image 
of the component $C_0^{\mathrm{sm}}$ of $X$. By \cite{EGAIV}, 
Proposition 2.8.1, $X_n$ (which is a scheme-theoretic closure) is flat over 
$\Spec ({\calO})$. The assumption 3 insures that $X_n$ is smooth over 
$K$. To prove the smoothness of $X_n$ over $\cal O$ it is therefore 
enough to prove that its special fiber is smooth over $k$, and this in turn 
follows from the assumption 4.
\end{proof}
\section{Bad reduction of Shimura curves.}
In this section we specialize the above considerations to the case of
Shimura curves. We begin by recalling some basic facts. Let $D$ be a product 
of an even number of distinct (finite) primes, let $N$ be a integer prime to $D$, 
and ${\cal O}^D_0 (N)$ an order of level $N$ in the quaternion algebra $B_D$ 
over $\QQ$ of discriminant $D$. As $B_D$ is indefinite (i.e. $\RR\otimes B_D$ 
is split), there is an embedding of the subgroup of ${{\cal O}_0^D (N)
}^\times$ consisting of elements with norm 1 into ${\mathrm{SL}}_2 (\RR )$, 
and the quotient of the Poincar\'e upper half-plane by this subgroup via the 
usual homographic action is a compact Riemann surface called the Shimura 
curve associated to ${\cal O}_0^D (N)$. We denote this curve by $X_0^D 
(N)(\CC )$. It follows from the work of Shimura that it has a 
canonical model over $\QQ$, and thanks to Deligne-Rapoport and Buzzard 
(model at primes dividing the level $N$), Cherednik-Drinfeld et al. (model at 
primes dividing the discriminant $D$) one has a model for $X_0^D (N)$ over
$\ZZ$. This model is constructed as the coarse moduli scheme
parametrizing abelian surfaces with multiplication by ${\cal O}_0^D 
(N)$, i.e. abelian varieties of dimension 2 with multiplication by a maximal
order of $B_D$ and a compatible $N$-level structure (plus an
additional technical assumption in characteristics dividing $D$ that we
will not discuss here). We will only need the description of
Cherednik-Drinfeld, i.e. we will concentrate on the localization of
$X_0^D (N)$ at primes dividing $D$. Fix such a prime $p$. To simplify
the notations, set $X:=X_0^D (N)$. According to Cherednik and Drinfeld,
$X$ is an admissible curve over $\ZZ_p$. Let ${\cal G} :={\cal G}
(X_{\FF_p} )$ be the dual graph of $X_{\FF_p}$ as in section 2. This 
graph is easier to describe when $D$ is a product of two primes $p$ and 
$q$, and from now on we assume that $D=pq$. Note however that dealing 
with the general case should not present conceptual difficulties.

In our case $D=pq$ it follows from the work of Ribet (\cite{Ribet},
Proposition 4.4 and 4.7) that the set $\calV$ of vertices of $\calG$
consists of two copies $\calV_1$ and $\calV_2$ of the set of
isomorphism classes of pairs $(E,C_N)$ where $E$ is a supersingular
elliptic curve over $\ol\FF_q$ and $C_N \subset E$ is a cyclic
subgroup of order $N$ (notice the switch between $p$ and $q$). The set 
$\calE$ of edges of the graph consists of $p$-isogenies between such 
objects. In other words the elements of $\cal E$ correspond to the 
isomorphism classes of triples $(E,C_N,C_p)$ where $C_p$ is a cyclic 
subgroup of order $p$, with source $(E,C_N )$ in $\calV_1$ and 
target $(E/C_p,C_N+C_p/C_p )$ in $\calV_2$. Note that this fixes an 
orientation on ${\cal G}$. The length of an edge corresponding to 
$(E,C_N,C_p)$ is ${\mathrm{card}} ({\mathrm{Aut}} (E,C_N,C_p )
/\pm 1)$.

For each $d$ dividing $DN$ such that ${\mathrm{gcd}}(d,DN/d)=1$, there 
is a distinguished involution $w_d$ on the curve $X_0^D (N)_{\ZZ_p}$,
a so-called Atkin-Lehner involution. The involution $w_q$ acts on 
$\cal G$ by sending $(E,C_N,C_p)$ to
$\Frob_q (E,C_N,C_p)$ so if $C_q$ is the kernel in $E$ of $\Frob_q$
then $\Frob_q (E,C_N,C_p)=(E/C_q,C_N +C_q/C_q, C_p+C_q/C_q)$. Note
that $w_q$ preserves the bipartition of $\calG$. The involution $w_p$
maps an edge $(E,C_N,C_p)$ to $-(E/C_p ,C_N+C_p /C_p ,E[p]/C_p )$,
therefore $w_p$ switches the bipartition of $\cal G$. Any quotient
$X^+ :=X/w_d$ is also an admissible curve and ${\cal G}^+ :={\cal G}
(X_{\FF_p} /w_d )={\cal G}(X_{\FF_p} )/w_d$ (taking into account the
ramification in the length of the edges).

Let $J:=\mathrm{Jac} (X)$ and $J^+ :=\mathrm{Jac} (X^+ )$ be the
N\'eron models over $\Spec (\ZZ_p )$ of the jacobians of $X_{\QQ_p}$
and $X^+_{\QQ_p}$ respectively. It follows from the previous section
that the neutral components $J_{\FF_p}^0$ and $J^{+0}_{\FF_p}$ of
$J_{\FF_p}$ and $J^+_{\FF_p}$ are tori, and the character group of
$J_{\FF_p}^0$ (respectively, $J^{+0}_{\FF_p}$) is isomorphic to $H_1
(\calG,\ZZ)$ (respectively, $H_1(\calG^+ ,\ZZ)$).

Let $\TT$ be the Hecke algebra of $\Gamma_0 (pqN)$. The results of 
Jacquet-Langlands imply that the elements of $\TT$ define correspondences 
on the Shimura curve $X=X^{pq}_0 (N)$. In order to formulate a version of 
Mestre-Oesterl\'e's graph method for Shimura curves, we need to understand 
the induced action of $\TT$ on the special fiber of $J$. To achieve this we 
use a kind of ``explicit twisted Jacquet-Langlands'' which is at the heart of 
\cite{Ribet}. In this article Ribet compares the character group of 
$\mathrm{Jac} (X_0^{pq} (N))_{\FF_p}$ with a subgroup of the character 
group of $\mathrm{Jac} (X_0 (pqN))_{\FF_q}$. To explain this result we 
introduce a few new notations. Let $\cal X$ be the group of degree $0$ 
divisors on the set of supersingular points in $X_0 (Nq)(\FF_{q^2} )$ and 
let $\cal L$ be the group defined in the analogous way for $X_0 (Npq)$. 
The two natural degeneracy maps $\alpha$ and $\beta$ from $X_0 (Npq)$ 
to  $X_0 (Nq)$ induce maps $\alpha_*$ and $\beta_*$ from $\cal L$ to 
$\cal X$. Call $Y$ the kernel of the morphism $(\alpha_* \oplus 
\beta_* )\colon {\cal L}\to {\cal X}\oplus {\cal X}$. Ribet shows 
that $\TT$ restricts to an action on $Y$ which cuts-out the $pq$-new quotient 
of $\TT$ (\cite{Ribet}, Theorem 3.20). Then he proves the following crucial 
result:
\begin{thm}
{\bf (Ribet)} There is a $\TT$-isomorphism (which is compatible with
the monodromy pairings) between $Y$ and the character group of
$\mathrm{Jac} (X^{pq}_0 (N) )^0_{\FF_p}$.
\end{thm}
\begin{proof}
This is \cite{Ribet}, Theorem 4.1.
\end{proof}
Note that the fact that $Y$ is the character group of the jacobian of a Shimura 
curve already followed implicitly from the description of the graph given above.
\subsection{Lemmata on component groups.}
We close this section by proving technical lemmas on the asymptotic shape and 
the order of components of a Shimura curve $X^{pq}$ (or its quotient $X^{pq} 
/w_q$) embedded in its jacobian. We fix some more notations. Let 
$S=\{ C_0 ,\dots ,C_g \}$ be the set of supersingular $j$-invariants 
in characteristic $q$. Write the bipartite graph ${\cal G}(X^{pq}_{\FF_p} )$
as a union ${\cal V}_1 \cup {\cal V}_2$. According to Ribet's description 
of the previous section, one may enumerate the vertices of ${\cal V}_r$ as 
$C_{r,0} ,\dots ,C_{r,g}$ for $r=1,2$. Each $C_{s} \in S$ corresponds to an
elliptic curve $E_{s}$, and we set $w(C_{s} ):=\mathrm{card} (\mathrm{Aut} 
(E_s )/(\pm 1))$. One knows that two $w(C_i )$ at most are different from 1, in 
which case they are equal to 2 or 3. Recall that the Eisenstein vector for $X_0  
(q)_{\overline{\FF}_q}$ is ${\mathrm{Eis}} :=(w(C_0 )^{-1}, w(C_1 )^{-1}, 
\dots ,w(C_g )^{-1})^t \in \frac{1}{12} \ZZ^S$. For $v=(v_{i} )_{i\in S} 
\in \CC^S$ we define the weight $w(v)$ as $\sum_{i\in S} v_{i}$. The weight 
$w({\mathrm{Eis}})=\sum w(C_r )^{-1}$ is roughly $g(X_0 (q)) \simeq q/12$.  
\begin{lem}
Fix a prime $q>3$. As the prime $p$ tends to infinity, the edges of
${\cal G}(X^{pq}_{\FF_p} )$ are equidistributed in the following sense. 
If $C_{1,j}$ and $C_{2,i}$ are elements of ${\cal V}_1$ and ${\cal V}_2$ 
respectively, the number of edges from $C_{1,j}$ to $C_{2,i}$ is
$$\frac{(p+1)}{w({\mathrm{Eis}})} \times \frac{1}{w(C_{1,j} )w
(C_{2,i} )} +O_q (\sqrt{p} ).$$
\end{lem}
\begin{proof}
We first deal with the generic case: when $w(C_{1,j} )=1$. Considering
the transpose $T_p \in M_{g+1} (\ZZ )$ of the $p^{\mathrm{th}}$-Brandt 
matrix for $\Gamma_{0} (q)$ with respect to the basis associated to $S$, one sees
that the number of edges from $C_{1,j}$ to $C_{2,j}$ is the $i^{\mathrm{th}}$
coefficient of the $j^{\mathrm{th}}$ column vector of $T_p$ (cf. for instance 
Proposition 4.4 of \cite{Gross}). Write $C_j =(0,\dots ,0,1,0,\dots ,0)^t$ in 
$\ZZ^S$, where the coefficient $1$ occurs at the $j^{\mathrm{th}}$ place. 
Decompose $C_j$ as $(1/w ({\mathrm{Eis}})) {\mathrm{Eis}} +c_j$, with 
$c_j$ a vector of weight 0. Now the lemma follows from the fact that $T_p 
({\mathrm{Eis}} )=(p+1){\mathrm{Eis}}$, together with Weil's bounds on 
cuspidal eigenvalues (i.e. if $a_{p}$ is such an eigenvalue, with $p$ prime to 
the level, then $| a_{p} |\leq 2\sqrt{p}$, see e.g. \cite{DiamondIm}, 
5.0.1 and references therein).

When $w(C_{1,j} )=2$ or $3$, one uses the same arguments, plus
the fact that, on a supersingular elliptic curve in characteristic $q$ 
corresponding to $C_{1,j}$, the $p$-isogenies which are conjugated by the 
exceptional automorphism $\zeta$ (=$\zeta_4$ or $\zeta_3$) are 
isomorphic. Identifying the set of $p$-isogenies with $\PP^1 (\FF_{p})$ 
one sees that, depending on the decomposition of $p$ in $\QQ (\zeta )$, 
the exceptional isomorphism may fix zero, one or two of these isogenies; 
the others belong to orbits of length $w(C_{1,j} )$. However the quantitative 
form of the statement does not depend on these different cases.
\end{proof}
\begin{lem}
Assume that $p$ and $q$verify the hypotheses of Theorem 1.1, case 2 with 
$p, q>3$. Let $C$ be a vertex of ${\cal G}((X^{pq} /w_q )_{\FF_p} )$ 
associated to the $\mathrm{Gal} (\FF_{q^2} /\FF_q )$-orbit of a 
supersingular invariant $j$. The number of edges emanating from $C$ is:
\begin{enumerate}
\item $(p+1)$ if $j$ belongs to $\FF_{q^2} \backslash \FF_q$;
\item $(p+1)/2$ if $j\in \FF_q$ and $j\not\equiv 0, 1728$;
\item $(p+3)/4$ if $j\equiv 1728$;
\item $(p+3\pm 2)/6$ if $j\equiv 0$.
\end{enumerate}
Setting $\varepsilon (C_{1,s} ,C_{2,r} )=1$ if $C_{1,s}$ or $C_{2,r}$
belongs to $\FF_{q^2} \backslash \FF_q$, and $\varepsilon (C_{1,s}
,C_{2,r} )=2$ if both lie in $\FF_q$, the number $N_{C_{1,s},C_{2,r}}$ 
of edges between two vertices $C_{1,s}$ and $C_{2,r}$ belonging to two 
different sets of the bipartition of ${\cal G}((X^{pq} /w_q )_{\FF_p} )$ 
is:
$$N_{C_{1,s},C_{2,r}} =\frac{p+1}{w({\mathrm{Eis}})} \times 
\frac{1}{\varepsilon (C_{1,s} ,C_{2,r} )w(C_{1,s} )w(C_{2,r} )} +
O_q (\sqrt{p} ) .$$
\end{lem}
\begin{proof}
The first part of the lemma follows form Ribet's description of the graph 
(sketched before Theorem 3.1). Indeed, let ${\cal C}$ be a vertex of 
${\cal G}(X^{pq}_{\FF_p} )$ associated to a supersingular curve $E$ with
invariant $j$. We recall that the involution $w_{q}$ acts as ${\mathrm{Frob
}}_{q}$ on ${\cal G}(X^{pq}_{\FF_p} )$. If $j$ belongs to $\FF_{q^2} 
\backslash \FF_q$, there are $p+1$ edges emanating from ${\cal C}$ in 
${\cal G}(X^{pq}_{\FF_p} )$, and as $w_{q} ({\cal C}) \neq {\cal C}$, 
the same is true for its image $C$ in ${\cal G} ((X^{pq} /w_q )_{\FF_p} )$, 
so case 1 of the lemma follows. Suppose that $j$ belongs to $\FF_q$, in 
which case $\cal C$ is fixed by $w_{q}$. The hypotheses on $p$ and $q$ 
insure that that no edge of ${\cal G}(X^{pq}_{\FF_p} )$ is fixed under the 
action of $w_q$ (see \cite{RotSkoYaf}, proof of Theorem 3.1 and 
\cite{OggMR}, p. 204 et seq.). This allows to conclude in case 2 of Lemma 
3.4. If $j\equiv 1728$ (respectively, $j\equiv 0$), one must take into 
account the fact that $E$ has an extra automorphism of order 4 (respect. of 
order 6) as in the proof of Lemma 3.3. As we assume that $p\equiv 1 \mod 4$, 
this shows that there are only $(p+3)/2$ (respect. $(p+3\pm 2)/3$) edges 
emanating from ${\cal C}$, and half this number from $C$. This proves the 
cases 3 and 4 of Lemma 3.4. The assertion about $N_{C_{1,s},C_{2,r}}$ is 
clear from Lemma 3.3. 
\end{proof}
To prove the next lemma we need an explicit description of component
groups. To achieve that purpose we recall the method of the proof of
\cite{BasRamanujan}, Proposition 9.2. See also \cite{BasBertDar},
section 1. 
 
  Let $\cal G$ be a graph with the set of vertices and edges $S_{0}$ and
$S_{1}$ respectively, and write $\ZZ^{S_{0}} [+]$ for the kernel of 
the ``weight'' map $f \mapsto \sum_{v\in S_{0}} f(v)$ from 
$\ZZ^{S_{0}}$ to $\ZZ$. Raynaud's results show that the component 
group $\Phi$ of the jacobian of an admissible curve with the graph $\cal G$ 
is the cokernel of the map $\iota \colon \ZZ^{S_{0}} \to \ZZ^{S_{0}} 
[+]$ given by the intersection matrix of $\cal G$. More precisely, if 
$\Delta \colon \ZZ \to \ZZ^{S_{0}}$ is the diagonal map, one has the 
exact sequence:  
$$(R){\ \ \ \ \ } 0\to \ZZ \stackrel{\Delta}{\to} \ZZ^{S_{0}} 
\stackrel{\iota}{\to} \ZZ^{S_{0}} [+] \to \Phi \to 0$$
(see e.g. \cite{BasRamanujan}, loc. cit.). Now having fixed an 
orientation on $\cal G$, let $s$ and $t$ be the  obvious ``source'' and 
``target'' maps from $S_{1}$ to $S_{0}$. Define $t^* \colon \ZZ^{
S_{0}} \to \ZZ^{S_{1}}$ by $(t^* f)(e)=f(t(e))$ and $t_{*} \colon 
\ZZ^{S_{1}} \to \ZZ^{S_{0}}$ by $(t_{*} f)(v)=\sum_{t(e)=v} f(e)$. 
Define similarly $s^*$ and $s_{*}$ and set $d^* :=t^* -s^*$ and $d_{*} 
:=t_{*} -s_{*}$. One readily checks that the map $\iota$ can be 
factorized as $\iota =-d_{*} d^*$ (cf. for instance 
\cite{DavSarnVal}, section 1.2, proof of Theorem 1.2.3). 
  
    The basic problem we discuss in the next lemma is the following: 
given an element in $\ZZ^{S_{0}} [+]$ of the form $(p+1)(C_{1} -
C_{2})$, is it in the image of $\iota$, or not? Recall that for $X$ an 
admissible curve, $\tilde{X}$ denotes the regular curve obtained from 
$X$ by blowing-up. 
\begin{lem}
Fix a prime $q$ such that $g(X_0 (q))\geq 5$ and suppose that $p \equiv 5
\mod 12$ and $q\equiv 3 \mod 4$. For $p$ large enough (compared to $q$), 
if a non-trivial difference between two irreducible components of
 $(\widetilde{X^{pq} /w_q})_{\FF_p}$ is killed by $(p+1)$ in the 
 component group of $\mathrm{Jac} (X^{pq} /w_q )_{\FF_p}$, then 
 these two components belong to the same set ${\cal V}_1$ or ${\cal 
 V}_2$ of the bipartition of the set of  components. 

  In particular, if $\cal J$ denotes the exceptional component arising 
from multiplication by $\zeta_4$ and $J$ is any other component, 
then $({\cal J}-J)$ is not killed by $p+1$.
\end{lem}
\begin{proof}
We consider the graph ${\cal G}^+ :={\cal G}((\widetilde{X^{pq} /
w_q} )_{\FF_p})$, whose sets of vertices and edges are still denoted by 
$S_{0}$ and $S_{1}$ respectively. As we assume that  $p\equiv 5\mod 
12$ and $q \equiv 3\mod 4$, we know that there is exactly one exceptional 
component $\cal J$ arising from multiplication by $\zeta_4$, and no 
exceptional component arising from multiplication by  $\zeta_3$ (see 
\cite{OggMR}, pp. 205 et seq., or the discussion preceding Theorem 4.1 
and the proof of Lemma 5.1 of the present paper. Note that this does not 
prevent a {\it vertex} from having ``multiplication by $\zeta_{3}$", 
i.e. from corresponding to a supersingular invariant $j\equiv 0 \mod q$). 
Let $w:=w({\mathrm{Eis}})$ be the weight of the Eisenstein vector for 
$X_0 (q)_{\FF_q}$ as in the Lemma 3.3. The hypothesis $g(X_0 (q))\geq 
5$ implies $w>4$. For $C\in S$ write $\varepsilon (C) =2$ if $C\in 
\FF_q$, and $\varepsilon (C)=1$ otherwise. Recall that the set of vertices 
of ${\cal G}^+$ other than $\cal J$ is ${\cal V}_1 \cup {\cal V}_2$. 
Call $J_1$ and $J_2$ the components surrounding $\cal J$. One has $w(
J_{1})=w(J_{2})=2$. There are exactly two edges emanating from ${\cal 
J}$, one towards $J_{1}$ and one towards $J_{2}$. Let $J$ be a component 
(possibly $J_1$) belonging to the set ${\cal V}_1$, say. We divide the set 
of vertices of ${\cal G}^+$ 
other than $\cal J$, $J$, $J_1$, $J_2$ into the subsets $\{ j_{i,1} \}$ and 
$\{ j_{i,2} \}$ of ${\cal V}_1$ and ${\cal V}_2$ respectively. By 
Lemma 3.4 if $(C_{1},C_{2})$ belongs to ${\cal V}_1 \times {\cal 
V}_2$ then $N_{C_{1},C_{2}} >0$ for $p\gg q$. We suppose this to be the 
case in what follows. (At this point, it may help the reader to make a rough 
sketch of what ${\cal G}^+$ looks like.) The symbol $O_q (\sqrt{p} )$ 
below will denote a priori different functions, which depend on the vertices of 
our graph, but as $q$ is fixed the set of vertices is fixed too, therefore those 
functions can be bounded by a single one. To prove Lemma 3.5 we have 
several kinds of generic configurations to consider for our differences $(C_1 -
C_2 )$, relative to the above data. Using the previous lemmas the computations 
are elementary but lengthy so we treat in detail only the cases of the last 
statement of the lemma and we leave other cases to the reader. (Note that in 
the computations below, we only need $w>3$, but the assumption $w>4$ in 
our statement is necessary for the cases we do not write down here). From now 
on we assume that there is a $v\in \ZZ^{S_{0}}$ such that $d_{*} d^{*} (v)
=(p+1)(J -{\cal J})=:i$. The discussion before the lemma shows that we need 
to obtain a contradiction for $p\gg q$.  

   We point out that our problem has a simple electrodynamical 
interpretation. Indeed, consider ${\cal G}^+$ as an electric circuit in 
which every edge has a resistance of one (ohm, say). The element $i\in 
\ZZ^{S_{0}} [+]$ may be seen as a vector of current: at each node 
$C\in S_{0}$, a current of $i(C)$ amp\`eres enters the circuit. If
$v\in \QQ^{S_{0}}$ satisfies $d_{*} d^{*} (v)=i$ then $v$ can be 
interpreted as a distribution of potential in volts at each node, 
corresponding to the current $i$ on ${\cal G}^+$. From this point of 
view our problem boils down to proving that no {\it integer valued} 
distribution of potential on the nodes allows a current of $p+1$ amp\`eres 
to enter the circuit at $J$ and leave it at $\cal J$. This electrodynamical 
interpretation is a common practice among graph theorists (see for instance 
\cite{Bollobas}, II.1). It is also used in \cite{BasRamanujan}, 
Proposition 9.2.  
\begin{sublem} 
{\bf (``law (K)").} With the above assumptions and notations, at each 
vertex $C$ in $S_{0}$ one has $\sum_{D \mapsto C} (v(C)-v(D))=
i(C)$, where ``$\sum_{D \mapsto C}$" means the sum over the 
neighbors of $C$ in ${\cal G}^+$.
\end{sublem}
{\bf Proof of 3.5.1.}
Clear from the fact that $i=d_{*} d^{*} (v)$.
$\square$

${ }$

   (Note that, from the electrodynamical point of view, this is just the 
conjunction of Ohm's law and Kirchhof's law (see \cite{Bollobas}, loc. 
cit.). Indeed, Ohm's law says that if a current $c$ passes through a wire 
of resistance $r$ between two nodes $A$ and $B$ with a difference in 
potential equal to $u$, then $u=r.c$. On the other hand, Kirchhoff's 
law says that the sum of currents passing through a node is zero. In our
case each wire in $S_{1}$ has resistance $r=1$ and $d^* (v)\in \ZZ^{
S_{1}}$ is the distribution of current through the wires of 
${\cal G}^+$.)  
\begin{sublem}
With the above hypotheses and notations, one has $v({\cal J})<v(C)<
v(J)$ for every vertex $C$ different from $J$ and ${\cal J}$. 
\end{sublem}
{\bf Proof of 3.5.2.} Let $C_{0}$ be such that $v(C_{0})$ is maximal 
among the $v(C)$ with $C\neq J$ and suppose that $v(C_{0})\geq v(J)$. 
By connexity of ${\cal G}^+$ one may assume that the vertex $C_{0}$ 
has a neighbor $D$ such that $v(D) <v(C_{0})$ (otherwise $v$ would be 
constant on ${\cal G}^+$, a contradiction with the exact sequence $(R)$ 
before Lemma 3.5). Therefore $\sum_{D\mapsto C_{0}} (v(C_{0})-v(
D))>0$, which is a contradiction with (3.5.1) as $i(C_{0})\leq 0$. The 
other inequality of the sublemma is obtained by symmetry. $\square$
  
${ }$

  (Note that from the electrodynamical point of view, the sublemma just
says that the distribution of potential must decrease between the nodes
$J$ (where the current enters) and $\cal J$ (where it comes out)). 

   We finally turn to the proof of the Lemma 3.5 itself, examining the 
different configurations that may occur for $J$ and $\cal J$. As we feel that 
the electrodynamical language gives slightly shorter and clearer proofs, we 
will use it freely. For readers who prefer more formal proofs however, we 
treat our first case below in both abstract and physical languages. It is our 
hope that this will clarify our method for the rest of the computations.   

Let ${\cal J}_2 \in \{ J_{2} \} \cup {\cal V}_{2}$ be such that $v(
{\cal J}_2 )=\max (v(j_{i,2} ),v(J_{2} ))$. Lemma 3.4 insures that there 
are $(p+1)/\varepsilon (J)w(J)$ edges leaving $J$. Recall that we assume 
that for each pair of vertices $(C_{1},C_{2})$ in ${\cal V}_1 \times 
{\cal V}_{2}$ the number $N_{C_{1},C_{2}}$ of edges between $C_1$ 
and $C_{2}$ is nonzero. By (3.5.2) one has $v(J)>v(J_{2} )$, $v(J)>v(j_{
i,2} )$ for all $j_{i,2}$ in ${\cal V}_{2}$, and 
by (3.5.1), $(p+1)=N_{J,J_{2}} (v(J)-v(J_{2})) +\sum_{j_{i,2}} N_{J,
j_{i,2}} (v(J)-v(j_{i,2}))\geq (v(J)-v({\cal J }_{2})) (N_{J,J_{2}} +
\sum_{j_{i,2}} N_{J,j_{i,2}} )\\
=(v(J)-v({\cal J}_{2})) (p+1)/\varepsilon (J)w(J)$, therefore $v(J)-v
({\cal J}_2 ) \leq \varepsilon (J)w(J)$.  

{\bf Case 1: Assume that $J\neq J_1$.}

{\bf 1a:  Suppose that $v$ is constant on $\{ J_{2} \} \cup {\cal V
}_{2}.$}

{\it (Abstract language.)} 
One has $v(J)-v({\cal J}_2 )= v(J)-v(J_2 )=\varepsilon (J)w(J)$ by the 
preceding computation. Applying the law $(K)$ (3.5.1) at each $j_{i,1}$ in 
${\cal V }_{1}$ shows that the $v(j_{i,1})$ are constant, equal to $v(J_{2})$. 
Applying (3.5.1) again at $J_{1}$, then at each $j_{i,2}$, then at $J$, gives 
that $v(J_{1})-v({\cal J})=(v(J_{2})-v(J_{1}))N_{J_{2},J_{1}}+\sum_{
j_{i,2}} N_{j_{i,2},J_{1}} (v(j_{i,2}) -v(J_{1})) \geq \sum_{j_{i,2}} N_{
j_{i,2},J_{1}} (v(j_{i,2}) -v(J_{1})) \\
\geq \sum_{j_{i,2}} N_{J,j_{i,2}} (v(J)-v(j_{i,2}))=(p+1)-(v(J)-v(J_{2}))
N_{J,J_{2}} =(p+1)(1-\frac{\varepsilon (J)}{2\varepsilon (J,J_{2})w
} )+O_q (\sqrt{p} )$ by Lemma 3.4. Now $v(J_{2})-v({\cal J})>v(J_{1}) 
-v({\cal J})$ so $(3.5.1)$ applied at $J_{2}$ gives $(p+1)(1-\frac{
\varepsilon (J)}{2\varepsilon (J,J_{2})w} )+O_q (\sqrt{p} ) \leq 
(v(J)-v(J_{2}))N_{J_{2},J} =(p+1)\frac{\varepsilon (J)}{2\varepsilon 
(J,J_{2})w} +O_q (\sqrt{p} )$. This implies $w\leq 1$ for $p\gg q$.  

{\it (Electrodynamical language.)}  
By the law $(K)$ at each $j_{i,1}$ in ${\cal V}_{1}$, one sees that the 
$v(j_{i,1})$s are constant and equal to $v(J_{2})$. Therefore all the current 
arriving from $J$ to the $j_{i,2}$s goes to $J_1$. It then goes to $\cal J$,
and this total current is at least $(p+1)(1-\frac{\varepsilon (J)}{2
\varepsilon (J,J_{2})w} )+O_q (\sqrt{p} )$ by Lemma 3.4. This 
means that $v(J_1 )-v({\cal J})$ is larger than this latter quantity, and the 
same is a fortiori true of $v(J_2 )-v({\cal J})$. But the current that $J_2$
receives from $J$ is only $(p+1)\frac{\varepsilon (J) }{2\varepsilon 
(J,J_{2})w} +O_q (\sqrt{p} )$, and this implies $w\leq 1$ for $p\gg q$.

{\bf 1.b: Suppose that $v$ is not constant on $\{ J_{2} \} \cup {\cal V
}_{2}$.} 

We first note that we have $\varepsilon (J)=2$, because if $J$ is quadratic 
over $\FF_{q}$ there are $(p+1)$ edges emanating from $J$, therefore (3.5.1) 
implies that $v$ is constant and equal to $v(J)-{1}$ on $\{ J_{2} \} \cup 
{\cal V}_{2}$. As there is a $j_{2}$ in $\{ J_{2} \} \cup {\cal V}_{2}$ 
such that $v(j_{2})<v({\cal J}_2 )$ one sees that $v({\cal J}_2 )>v(J_{1})$ 
(otherwise we get a contradiction with the fact that $J_{1}$ has no current to 
offer to $j_{2}$), and similarly $v({\cal J}_2 )>v(j_{i,1})$ for all $j_{i,1}$ 
in ${\cal V}_{1}$. Therefore the node ${\cal J}_2$ has to deal at least 
$\frac{p+1}{\varepsilon ({\cal J}_2 )w({\cal J}_2 )} (1-\frac{1}{w.
w(J)}) +O_{q} (\sqrt{p} )$ amp\`eres, whereas it receives (from $J$) 
strictly less than $\frac{(p+1)}{\varepsilon ({\cal J}_2 )w({\cal J}_2 
)w(J)w} (\varepsilon (J)w(J)) +O_{q} (\sqrt{p} )$. This gives $w\leq 
\varepsilon (J) +1/w(J)\leq 3$ for $p\gg q$.

{\bf Case 2. Assume that $J=J_1$.} As in the case 1.b, the node ${\cal J}_2$ 
has to deal at least $\frac{p+1}{\varepsilon ({\cal J}_2 )w({\cal J}_2 
)} (1-\frac{1}{2w}) +O_{q} (\sqrt{p} )$ amp\`eres, whereas it receives 
from $J_{1}$ less than  $\frac{(p+1)}{2\varepsilon ({\cal J}_2 )w({\cal 
J}_2 )w} (v(J_{1})-v({\cal J}_2 )) +O_{q} (\sqrt{p} )\leq \frac{2(p+1)
}{\varepsilon ({\cal J}_2 )w({\cal J}_2 )w} +O_{q} (\sqrt{p} )$. This 
gives $w\leq 5/2$ for $p\gg q$.
 
This completes the examination of the cases of differences of
components of shape ${\cal J}-J$, which will be the only ones of
interest to us in the following of the paper. The other cases of the
lemma are similarly elementary and lengthy to treat.
\end{proof}
\begin{rem} The above makes it clear that, if $p$ is large enough
relatively to $q$, the graphs ${\cal G}$ or ${\cal G}^+$ are
``non-disconnecting''. It follows that the maps from $X^{pq, {\mathrm{
sm}}}_{\ZZ_p}$ and $(X^{pq}/w_q )^{\mathrm{sm}}_{\ZZ_p}$ to their 
jacobian over the same basis are closed immersions, by Proposition 9.5 of 
\cite{BasRamanujan}. We will use this remark before the Theorem 5.3. 
\end{rem}
\section{Shimura curves and results of Gross, Zhang and Waldspurger.}
In this section we show how the results of Gross, Zhang and Waldspurger 
on the special values of automorphic $L$-functions describe the winding
quotient of the jacobian of our Shimura curves, generalizing 
\cite{ParentX0+}, Proposition 4.2.

Let $G_{\QQ}$ denote the absolute Galois group of $\QQ$ and 
$S_{2}^{\mathrm{new}} (\Gamma_{0} (N))_{\overline{\QQ}}$ the set 
of newforms over $\overline{\QQ}$ of weight 2 for $\Gamma_{0} (N)$. 
For any integer $N$, the theory of Shimura associates to the decomposition 
into  $G_{\QQ}$-orbits of $S_2^{\mathrm{new}} (\Gamma_0 (N))_{
\overline{\QQ}}$ a similar decomposition over $\QQ$, up to isogeny, of 
$J_0 (N)$ into a product $\prod_{G_{\QQ} -{\mathrm{orbits}}} A_f$. 
According to Shimura, Carayol et al., one has
$$L(J_0 (N),s)=* \prod_{G_{\QQ} -{\mathrm{orbits}}} L(A_f ,s)=*
\prod_{G_{\QQ} -{\mathrm{orbits}}} \prod_{(f^\sigma \in G_{\QQ} 
\cdot f)} L(f^\sigma ,s) $$
(where the $*$s denote non-zero constants). The Hecke algebra $\TT :=
\TT_{\Gamma_0 (N)}$ is a free $\ZZ$-module of finite type. The above
decomposition of $J_0 (N)$ corresponds to idempotents in $\TT \otimes
\QQ$ and there is a unique ideal $I_e$ of $\TT$ such that $J_0 (N)/I_e 
J_0 (N)$ is isogenous to $\prod_{L(A_f ,1) \neq 0} A_f$ and $I_{e}$ is 
saturated (i.e. equal to the inverse image of $I_{e} \otimes \QQ$ under the 
map $\TT \to \TT \otimes
\QQ$). This quotient abelian variety is called the winding quotient
and we denote  
it by $J_e (N)$. A milestone theorem of Kolyvagin-Logachev (see 
\cite{Kolyvagin}), generalized by Kato in \cite{Kato}, implies that $J_e (N)
(\QQ )$ is finite (this is part of a special case of the conjecture of
Birch and  
Swinnerton-Dyer). The theory of Jacquet-Langlands establishes a correspondence 
between automorphic forms for $\GL_{2}$ on the one hand, and automorphic 
forms for the multiplicative group of nonsplit quaternion algebras, on the
other hand (see  
e.g. \cite{JacquetLanglands}). This implies that one can also make the 
above construction for any (quotient of) a general Shimura curve instead 
of the classical $X_0 (N)$.

We want to develop a graph method on our Shimura curves in a similar
way to \cite{ParentX0+}. The $\ZZ$-modules $\cal L$ and $Y$ of 
the end of section 3.1 can be interpreted as character groups or, as we 
are working with tori, as spaces of differential 1-forms on the
corresponding geometric objects. The use of special values formulae
for suitable Rankin-Selberg $L$-functions allows us to describe in
these spaces the subspaces of differentials pertaining to the winding
quotient. To be more precise we need to fix some notations and recall
some results from Gross' theory (a neat introduction to this circle of
ideas can be found in \cite{VatsalRS}). If $M$ is a $\ZZ$-module, 
define $\hat{M} :=M\otimes \hat{\ZZ}$. Let $B$ be the
quaternion algebra over $\QQ$ which is ramified precisely at $q$ and
$\infty$. Choose a Eichler order $R$ of level $p$ in $B$, and let $\{
R_1 :=R,... ,R_n \}$ be a set of Eichler orders in $B$ corresponding to
representatives for ${\hat{R}}^* \backslash {\hat{B}}^* /B^*$: to 
a double coset $g:=(g_2 ,g_3 , \dots ,g_l ,\dots )$ we associate the 
$B^*$-conjugation class of the Eichler order $B\cap g^{-1} \hat{R} g$. 
Recall that ${\hat{R}}^* \backslash {\hat{B}}^* /B^*$ is in 
one-to-one correspondence with the set $\cal S$ of singular points of
$X_0 (pq)_{\FF_q}$ (see \cite{Ribet}, 3.3 and 3.4). The order $R_i$ 
associated to a point $x=(E,C_p )$ is such that $R_i \simeq {\mathrm{ 
End}}_{\overline{\FF}_q} (E,C_p )$.

If $L$ is a quadratic number field, it embeds in $B$ if and only
if its localization at ramification primes for $B$ is a field, i.e.
$L$ is a quadratic imaginary field in which $q$ is inert or ramified.
Then, for an order ${\cal O}$ of $L$, a morphism of algebras $\sigma 
: L\hookrightarrow B$, and a Eichler order ${\cal R}$ as above, the 
pair $(\sigma ,{\cal R})$ is said to be an {\it optimal embedding} of 
$\cal O$ in ${\cal R}$ if $\sigma (L)\cap {\cal R}=\sigma ({\cal 
O} )$. Such an embedding exists if and only if $p$ splits or ramifies in 
$\cal O$. If $-D$ is a negative integer, let $h(-D)$ be the class number 
of the quadratic order ${\cal O}_{-D}$ with discriminant $-D$ (if it
exists), and let $h_i (-D)$ be the number of optimal embeddings of
${\cal O}_{-D}$ in $R_i$ modulo conjugation by $R_i^*$. We define 
in $\ZZ^{\cal S}$ the element:
$$e_D :=\sum_{i=1}^{n} \frac{h_i (-D)}{\vert {\mathrm{Aut}} 
(R_{i}^* /\pm 1)\vert } [R_i ].$$
We call those elements {\it Gross vectors}.

Now let $f$ be a newform of weight 2 for $\Gamma_0 (pq)$. If $-D$ is
a quadratic imaginary discriminant as above, we write $-D=-d.c^2$,
where $-d$ is a fundamental discriminant and $c$ a conductor. Let
$\varepsilon_D$ be the non-trivial quadratic Dirichlet character
associated to $\QQ (\sqrt{-D} )$, and $f\otimes \varepsilon_D$ the
twist of $f$ by $\varepsilon_D$. One knows that the natural action
of the Hecke operators as correspondences endows ${\cal L}$
with a structure of $\TT$-module (the Brandt Hecke module), such
that ${\cal L} \otimes \QQ$ is a free $\TT \otimes \QQ$-module of  rank 1 
(see \cite{Gross}, section 5, for the case of prime level, and for instance 
\cite{Emerton} for more general results). Let $({\cal L} \otimes 
\overline{\QQ} )^f$ be the $\TT_{\overline{\QQ}}$-eigenspace in 
${\cal L}\otimes \overline{\QQ}$ associated to $f$, and let $e_{f,D}$ 
be the component of $e_D$ on $({\cal L} \otimes \overline{\QQ} )^f$.
\begin{thm} {\bf (Gross, Zhang, Waldspurger)} With notations as above,
the number $L(f,1)L(f\otimes \varepsilon_D ,1)$ is non-zero if and 
only if $e_{f,D}$ is non-zero.
\end{thm}
\begin{proof} 
This follows from a result of Waldspurger, as explained, in a slightly different 
language, in section 1.4 of \cite{CornutVatsal}. Under the assumption that 
$(pq,D)=1$ and $(c,pqD)=1$, see also \cite{VatsalRS}, Theorem 6.4. 
Note that, when $f$ has prime level, an exact formula expressing $L(f,1)L(f
\otimes \varepsilon_D ,1)$ in terms of $e_{f,D}$ was proved by Gross in 
his seminal work \cite{Gross}, Corollary 11.6. This was generalized in 
\cite{Zhang}, Theorem 1.3.2, or \cite{Zhang2}, Theorem 7.1. It was this 
formula that we used in \cite{ParentX0+}. But even the newest versions of 
this exact Gross-Zhang formula require mildly restrictive hypotheses. In the rest 
of this work however, as in the arithmetic applications of \cite{CornutVatsal}, 
one only needs the nonvanishing statement of Theorem 4.1. \end{proof}
We now come back to the case of modular curves. It follows
from \cite{Ribet}, Proposition 3.1, that the group ${\cal L}$ of 
section 3.1 is isomorphic as a $\TT$-module to the character group of 
the maximal subtorus of $J_0 (pq)_{\FF_q}$.  
\begin{prop}
Let $\EE$ be the $\QQ$-vector subspace of ${\cal L} \otimes \QQ$
spanned by the orthogonal projections of Gross vectors
(relative to the monodromy pairing $\langle \cdot ,\cdot \rangle$, 
cf. section 2). Then $\EE ={\cal L} [I_e ]\otimes \QQ$.
\end{prop}
\begin{proof}
This is a straightforward variant of \cite{ParentX0+}, Proposition 
4.2.
\end{proof}
We obtain an analogous result in the case of Shimura curves.
\begin{cor}
The subspace $\EE \cap Y_\QQ$ of $Y_{\QQ} =Y\otimes \QQ$ is equal 
to the character group of the winding quotient of $J_0^{pq-\mathrm{
new}} (pq)$, i.e. the character group of the winding quotient of
$\mathrm{Jac} (X^{pq})^0_{\FF_p}$.
\end{cor}
\begin{proof}
Combine Proposition 4.2 with Theorem 3.1.
\end{proof}
The Corollary 4.3 says that homology classes of the closed paths in the
graph of $X^{pq}_{\FF_p}$ made of Gross vectors are characters
of the winding quotient. This gives a concrete method for obtaining 
those characters, which is a crucial point of our strategy.
\section{General method.}
We turn to our initial problem of proving the triviality of the
rational points on Atkin-Lehner quotients of Shimura curves. Set as before
$X=X^{pq}/w_q$. Suppose that there exists a point $P$ in $X(\QQ)$, and 
let $P_0 :=w_p (P)$. Our aim is to prove that $P=P_0$, which implies
that $P$ is a special (or CM) point (cf. e.g. the first chapter of
\cite{Clark} or \cite{AlsinaBayer}, Remark 6.5) - or leads to a 
contradiction if there are no special points in $X(\QQ )$. Consider the 
following diagram of $\QQ$-morphisms:
$$
\begin{array}{ccccl}
& & X_{p+1} & & \\
& \nearrow & & \searrow & \\
X & \stackrel{\phi}{\to} & J & \stackrel{[p+1]}{\longrightarrow}
& J \\
& & {\scriptstyle \alpha} \downarrow & & \downarrow {\scriptstyle 
\alpha} \\
& & J_e^{w_q =1} & \stackrel{[p+1]}{\longrightarrow} & J_e^{w_q =1}
\end{array}
$$
where $\phi (Q)=(Q)-(P_0 )$, $X_{p+1} =[p+1]\circ \phi (X)$ as at the
end of section 2, $\alpha$ is the projection of $J={\mathrm{Jac}} (X)$ to 
its winding quotient $J_e^{w_q =1}$ (isogenous to $J_e /(1-w_q )J_e$ for 
$J_e$ the winding quotient of ${\mathrm{Jac}} (X^{pq} )$). It is clear that 
everything commutes in this diagram. We define $\psi =[p+1]\circ \phi$ 
and $\Psi =\alpha \circ \psi$.

Let $\tilde{X}$ be the desingularization of $X$ over $\ZZ_{p}$ as in
section 2. As $\tilde{X}$ is regular, any point in $\tilde{X} (\ZZ_{p} )$ 
belongs to $\tilde{X}^{\mathrm{sm}} (\ZZ_{p} )$. For such a 
$Q$ in $\tilde{X} (\ZZ_p )$, let $C_Q$ be the vertex of $\calG 
(\tilde{X})$ corresponding to the specialization component of $Q$.
\begin{lem}
Assume that $p$ and $q$ satisfy the hypotheses of Theorem 1.1, case 2. Then 
$\tilde{X}_{\FF_p}$ has a unique $\FF_p$-rational component which is 
necessarily $C_{P_0}$. This component intersects the other components in two 
points exactly. The embedding $\phi$ identifies $C_{P_0}^{\mathrm{sm}}$ 
with a one-dimensional subtorus $T$ of $J$ which is isomorphic to the 
non-trivial twist of $\GG_{m}$ over $\FF_p$ (cf. Proposition 2.1).
\end{lem}
In particular, for any $Q$ in $\tilde{X} (\ZZ_p)$, one has
$C_Q =C_{w_p (Q)}=C_{P_0}$.
\begin{proof}
We use again the work of Ogg (\cite{OggMR}), Th\'eor\`eme
of p. 206. In loc. cit. indeed, as we already noticed, it is explained that the 
vertices and the edges of ${\cal G}(X^{pq}_{\FF_p} )$ are defined over 
$\FF_{p^2}$, and that the action of $\mathrm{Frob}_p$ on components and 
singular points of $X^{pq}_{\FF_p}$ corresponds to the action of $w_p$ on the
graph ${\cal G}(X^{pq}_{\FF_p})$. Recall that the set of vertices of 
${\cal G} (X^{pq}_{\FF_p})$ consists of two copies of the set of isomorphism 
classes of supersingular elliptic curves in characteristic $q$ and that $w_p$
sends a vertex from the first copy to the symmetrical vertex from
the second. It follows that no vertex of ${\cal G} (X^{pq}_{\FF_p})$
is fixed by $w_p$. Passing to the quotient curve one sees that the
$\FF_p$-rational components of $X_{\FF_p}$ are exceptional, and more
precisely are obtained by blowing-up the edges $y$ of ${\cal G}
(X_{\FF_p})$ of length two and such that $w_p(y) = -y$ (i.e. the
edges which are preserved by $w_{p}$, but whose source and target
are permuted).

As already discussed, those exceptional components arise in two situations, the 
``cas ramifi\'e'' and the ``cas non-ramifi\'e'' of Ogg. In the first case the lifting 
of $y$ to ${\cal G}(X^{pq}_{\FF_p})$ consists of one edge fixed by $w_q$. 
In the second case $y$ has two preimages in ${\cal G}(X^{pq}_{\FF_p})$ 
exchanged by $w_q$, see loc. cit., pp. 204 and 205. Note that the fact that $pq$ 
is odd implies that $X^{pq} (\QQ_p )$ itself is empty: this follows from the point 
$i)$ of the Theorem of p. 206 of \cite{OggMR}. The hypotheses of the present 
lemma, that $\left( \frac{q}{p} \right) =-1$ and the congruences satisfied 
by $p$ and $q$, imply that we are in the ``cas non-ramifi\'e''. Let $C_1$ and 
$C_2$ be the liftings to ${\cal G}(X^{pq}_{\FF_p})$ of a $\FF_p$-rational 
component $C_{P_0}$ in ${\cal G} (X_{\FF_p})$. Using Ribet's description 
(\cite{Ribet}, Section 3) we see that the $C_i$'s correspond to Eichler 
orders $\Theta$ of level $p$ in $B_{q\infty}$ which contain a fourth root 
of unity $\zeta_4$, so that the length of $C_{i}$ is ${\mathrm{card}} 
({\mathrm{Aut}} (\Theta /\pm 1))=2$. By Eichler mass formula (e.g. 
\cite{Gross}, 1.12), as $q\equiv 3\mod 4$ there is (up to conjugacy) 
exactly one maximal order in $B_{q\infty}$ containing $\zeta_4$ and, as 
$p\equiv 1\mod 4$, exactly two such orders of level $p$. (One can see this 
by using Ribet's description according to which $\Theta$ is the 
endomorphism ring of an ``enhanced elliptic curve'', i.e. a pair consisting 
of a supersingular elliptic curve in characteristic $q$ and a cyclic subgroup 
of order $p$.) Those two orders are switched by $w_{q}$, and this shows the 
existence and uniqueness of $C_{P_0}$. As $C_{P_0}$ comes from a blow-up, 
it intersects other components at two singular points only. To prove the last 
statement of the lemma we remark that Proposition 2.1 gives the desired 
identification between $C_{P_{0}}$ and a one-dimensional subtorus $T$. 
Using again that the Frobenius action on ${\cal G}(\tilde{X}_{\FF_p} )$ 
is given by $w_p$, and that $w_p$ switches the two singular points on 
$C_{P_0}$, we finally obtain that $\mathrm{card} (T(\FF_p ))=p+1$.
\end{proof}
\begin{lem}
If $n$ is a integer such that $\phi (X)(\CC )\cap J[n] (\CC )\neq\{
0\}$, then $n\geq (pq)/245$ provided that $p$ is at least $19$ and $q$
is at least $245$.
\end{lem}
\begin{proof}
We first recall that a projective and smooth curve $\cal C$ over a field $k$ 
is said to be $n$-gonal over $k$ if there exists a $k$-morphism of degree $n$ 
from $\cal C$ to $\PP^1$. The gonality of $\cal C$ over $k$ is the smallest 
$n$ for which $\cal C$ is $n$-gonal over $k$. In the situation of the lemma, 
if $\phi (X)(\CC )\cap J[n](\CC )$ contains a non-trivial point $Q$, then 
$n((Q)-(P_{0}))$ is the divisor of a rational function which provides a 
morphism of degree $n$ to the projective line, therefore $X$ is $n$-gonal over 
$\CC$. This means that the $\CC$-gonality of $X^{pq}$ is at most $2n$. But 
according to \cite{Abramovich}, Theorem 1.1, this gonality is at least 
$(21/200)(g(X^{pq}) -1)$. The genus of $X^{pq}$ is given by the following 
formula (\cite{Shim}, chapter 2; see also \cite{Clark}, Proposition 46 for 
an immediate formula (note that $pq$ is odd), or \cite{SkoYaf}, remark after 
Proposition 1.4):
$$
g(X^{pq}) = 1 - \frac{1}{4}(1-\left( \frac{-1}{p}\right))
(1-\left(\frac{-1}{q}\right))
- \frac{1}{3}(1-\left( \frac{-3}{p}\right))(1-\left( \frac{-3}{q}
\right))+ \frac{(p-1)(q-1)}{12}
$$
hence
$$
g(X^{pq}) - 1 \geq -\frac{7}{3} + \frac{(p-1)(q-1)}{12} .
$$
An elementary computation now shows that
$n\geq \frac{21}{400}(g(X^{pq}) - 1) \geq \frac{pq}{245}$
provided that $p$ and $q$ are as in the statement.
\end{proof}
The two preceding lemmas, together with Lemma 3.5 and Remark 3.6, 
show that for $q$ and $p$ large enough and verifying $\left( \frac{q}{p} 
\right) =-1$, $p\equiv 5 \mod 12$, $q\equiv 3 \mod 4$, all the four 
hypotheses of Proposition 2.2 are verified with $n=p+1$: we may therefore 
speak of the model $X_{p+1}$ of $X_{p+1,\QQ_p}$ over $\ZZ_p$ provided 
by Proposition 2.2.
\begin{thm}
Keep the above notations and assumptions. Assume that there exists a
closed path on ${\cal G}(X)$ which is made of Gross vectors and
contains, with prime-to-$p$ multiplicity, the edge corresponding to
$C_{P_0}$. Assume furthermore that $q>245$ and $p\gg q$. Then
$P=P_0$ and hence $P$ is a special point.
\end{thm}
\begin{proof}
We choose a basis of $H_1 ({\cal G}(X) ,\ZZ )$ as in the proof of
Proposition 2.1, with only one vector $V_{P_0}$ containing (with multiplicity 
one) the edge corresponding to $C_{P_0}$. The path $\gamma$ made of the 
Gross vectors of the statement has a component at $V_{P_0}$ which is nonzero 
mod $p$, and $\gamma$ corresponds to a character of the winding quotient by 
Corollary 4.3. We use our  basis of $H_1 ({\cal G}(X) ,\ZZ )$ as in the proof 
of Proposition 2.1 and pick a $\FF_p$-parameter $x$ on $C_{P_0}$. Let $a$ 
and $\overline{a}$ be the values taken by $x$ at the singular points of 
$C_{P_0}$, which are $\FF_{p^2}$-conjugate (Lemma 5.1), and set 
$\lambda :=(x(P_{0} ) -\overline{a})/(x({P_{0}} )-a)$. With these 
notations, a point $Q$ in $C_{P_0}$ is mapped by $\phi$ to $(\lambda 
\times \frac{x(Q) -a}{x(Q) -\overline{a}} ,1,\dots ,1)$ (Proposition 
2.1). This gives the identification via $\phi_{\FF_p}$ of $C_{P_0}$ with 
the twisted form of $\GG_m$ over $\FF_p$ called $T$ in Lemma 5.1. As 
the fiber of $J$ at $p$ is purely toric, its cotangent space is naturally identified 
with the tensor product of its character group with $\FF_{p}$. We call 
$\omega_{\gamma}$ 
the invariant differential on $J_{\FF_{p}}$ corresponding to $\gamma$ via this 
identification. One now computes that, if $\beta$ is the component of 
$\gamma$ at $V_{P_{0}}$, the pull-back of $\omega_{\gamma}$ to the
cotangent space of $X_{\FF_{p}}$ at $P_{0}$ is $\beta (1/(x(P_{0} )-a) 
-1/(x(P_{0} )-\overline{a})) dx$, which of course does not vanish. The 
above shows that $\mathrm{Cot} (\phi )$ is nonzero at $\phi (P_{0} )$, and 
this means that $\phi$ is a formal immersion at $P_{0}$. We briefly indicate 
how one can get the same for $\alpha \circ \phi$, with $\alpha$ defined
at the beginning of section 5 (see e.g. the proof of \cite{ParentX0+}, 
Proposition 3.2). Let $I_{e}^{w_{q} =1}$ be the Eisenstein ideal relative to 
$J$. The exact sequence of abelian varieties over $\QQ$:
$$0\to I_{e}^{w_{q} =1} .J_{/\QQ} \to J_{/\QQ} \to J_{e/\QQ}^{w_{q} 
=1} \to 0$$
induces an exact sequence of cotangent spaces at the origin, which are free 
$\ZZ_{p}$-modules of finite rank:
$$0\to {\mathrm{Cot}}  (J_{e/\ZZ_{p}}^{w_{q} =1} )\to {\mathrm{Cot}} 
(J_{/\ZZ_{p}})  \to {\mathrm{Cot}} (I_{e}^{w_{q} =1} .J_{/\ZZ_{p}} )
\to 0$$
(this follows from the semistability of $J$, see the references in 
\cite{ParentX0+}, loc. cit.). Using our previous remarks on ${\gamma}$
and the identification between cotangent spaces and character groups, one sees 
that $\omega_{\gamma}$ actually belongs to the subspace ${\mathrm{Cot}} (
J_{e/\FF_{p}}^{w_{q} =1} )$ of ${\mathrm{Cot}} (J_{/\FF_{p}} )$. 
Therefore the nonvanishing of 
$\mathrm{Cot} (\phi )(\omega_{\gamma}) =\mathrm{Cot}  (\alpha \circ 
\phi )(\omega_{\gamma})$ shows that $\alpha \circ \phi$ also is a formal 
immersion at $P_{0}$ (and actually at any point on the component 
$C_{P_0}^{\mathrm{sm}}$). As multiplication by an integer $n$ on an abelian 
variety induces multiplication by $n$ on its cotangent space (with respect to the 
module structure), one sees that the map $\Psi$ defined at the beginning of 
this section is a formal immersion at $P_0$, because $p+1$ is prime to $p$. 
Finally, as $\Psi$ factors through $X_{p+1}$, the map $X_{p+1} \to J_e^{
w_q =1}$ is again a formal immersion at $(p+1)\phi (P_0 )(\FF_p )$. 
Furthermore, Lemma 5.1 shows that $(p+1)\phi (P)(\FF_p ) =0(\FF_p )=
(p+1)\phi (P_0 )(\FF_p ) $ in $X_{p+1} (\FF_p)\hookrightarrow J(\FF_p)$.

Summing up: one has two points $(p+1)\phi (P_0 :=w_p (P))$ and $(p+1)\phi 
(P)$ in $X_{p+1}(\QQ)$, specializing to the same point at $p$, where the natural 
map to a winding quotient is a formal immersion. This allows us to apply Mazur's 
classical argument to the curve $X_{p+1}$ (see e.g. \cite{ParentX0+}, proof of
the Proposition 3.1 for a version akin to the present one). We therefore obtain that
$(p+1)((w_p (P))-(P))=0$ in $X_{p+1} (\QQ ) \into J(\QQ )$. But this means that 
either $P=w_p (P)$, in which case $P$ is a special point, or $X$ is $(p+1)$-gonal. 
Now Lemma 5.2 shows that the choice of $p$ and $q$ implies that the latter case is 
impossible.
\end{proof}
\section{An explicit infinite family.}
We finally illustrate our method by considering the family of curves
$X^{p\cdot 251} /w_{251}$, for $p$ large and running through the
infinite set of primes verifying the explicit congruences of Lemma 6.1
below, which ensures that the hypotheses of Theorem 5.3 are verified. 
Using Gross vectors, we will construct explicit closed paths on the 
graph $\calG :=\calG ((X^{p\cdot 251} /w_{251} )_{\FF_p })$. 
Call ${\cal S}_1$ and ${\cal S}_2$ two copies of the supersingular
invariants in characteristic $251$. The vertices of $\calG$ are the
${\mathrm{Gal}} (\FF_{251^2} /\FF_{251} )$-orbits in the ${\cal
S}_i$s. Denote by $\tilde{\cal G} :={\cal G} (\widetilde{X^{p
\cdot 251} /w_{251}} )$ the desingularized graph. Each ${\cal S}_i$ 
is given by the roots of the supersingular polynomial in characteristic 
251, which is:
$$j(j+29)(j+19)(j-64)(j-4)(j-24)(j-30)(j-35)(j-101)(j+112)(j+66)
(j+52)(j+44)\times$$
$$(j+38)(j^2 -60j-81)(j^2 -81j-68)(j^2 -105j+116)(j^2 +93j +91)$$
(see for instance \cite{AntwerpIV}). The second factor of this product 
is congruent to $(j-1728) \mod 251$. Let $E_{1728}$ be an elliptic curve
over $\FF_{251}$ with $j$-invariant $1728$. As $p\equiv 1 \mod 4$, 
there are two embeddings $\ZZ [\zeta_{4} ]\hookrightarrow 
{\mathrm{End}} (E_{1728} , C_{p} )$, for two $p$-isogenies $C_{p}$.
They give rise to the two only edges of length 2 in $\calG (X^{p\cdot 
251} )$ from the vertex $j=1728$ in ${\cal S}_1$ to the vertex $j=1728$ 
in ${\cal S}_2$ (see the proof of Lemma 5.1). The Gross vector $e_{4}$ 
in the graph of $X^{p\cdot 251}$ corresponds to the sums of these two 
edges, which are switched by $w_{251}$, hence $e_{4}$ gives a single 
edge of length 2 in ${\calG} (X^{p\cdot 251} /w_{251} )_{\FF_{p}}$. 
In $\tilde{\cal G}$, this edge finally splits into the two edges (of length 
1) surrounding the exceptional vertex which corresponds to the only 
$\FF_p$-rational component of $(\widetilde{X^{p\cdot 251} /
w_{251}} )_{\FF_p }$ (see Lemma 5.1).

Now after looking in tables (e.g. \cite{Lario}) among the Heegner
points of degree one and two, one considers the following other three 
Gross vectors: $e_{28}$, $e_{36}$, $e_{267}$, associated to the orders 
$\ZZ [\sqrt{-7}]$, $\ZZ [3\sqrt{-1}]$ and $\ZZ  [(1+\sqrt{-267})/2]$ 
respectively, with class polynomials:
$$
P_{-28} =(j-16581375)\equiv (j-64) \mod 251 ;
$$
$$
P_{-36} =(j^2 -153542016*j -1790957481984)\equiv (j-64)(j+19) \mod
251 ;
$$
$$
P_{-267} =(j^2 +19683091854079488000000*j
+531429662672621376897024000000) $$
$$
\equiv (j+19)(j+29) \mod 251 .
$$
One readily checks that the corresponding quadratic orders embed into 
the quaternion algebra $B_{251\cdot \infty}$ (i.e. that $251$ is inert 
in these orders), and if $p$ splits in these orders, then for each one we 
obtain two $p$-isogenies linking the (one or two) $j$s of each 
polynomial. Those isogenies are interpreted as edges from ${\cal S}_1$ 
to ${\cal S}_2$ on the graph ${\cal G}(X^{p\cdot 251} )_{\FF_p}$, 
then on $\tilde{\cal G}$, in the same way as explained in the above 
case of $e_{4}$. We sum-up the above observations in the following:
\begin{lem}
Assume that $p$ verifies the following congruence conditions:
$$p\equiv 5 \mod 12, \left( \frac{-7}{p}\right)=1, \left(
\frac{-267}{p}\right)=1, {and}\ \left( \frac{p}{251}\right)=-1.$$
Then the sum of Heegner points: $(e_{4} -e_{267} +e_{36}
-e_{28} )$ defines a closed path in $\calG ((X^{p\cdot 251} /w_{251}
)_{\FF_p })$.
\end{lem}
\begin{proof}
A straightforward application of Theorem 3.1.
\end{proof}
\begin{rem}
The path of Lemma 6.1 may have different shapes depending on $p$ - for 
instance it can have one, two or three connected components. These shapes 
depend on whether or not the degree-two Heegner points link their two
different roots $j$. This in turn is determined by the class of the
ideals above $p$ in the Picard group of each corresponding quadratic order, 
by the theory of complex multiplication for elliptic curves. In any case 
however, one gets a closed path passing through the edge having CM by 
$\zeta_{4}$, which is the only relevant point regarding our problem.
\end{rem}

   Consider the family of curves $X^{p \cdot 251} /w_{251}$ with $p$ 
satisfying the congruences above. It follows from Theorem 5.3 that the only 
possible $\QQ$-rational points on the curves of this family are special points.

Actually, the curves in our family do {\it not} have any special rational 
point either. For a curve $X^{p \cdot 251}/w_{251}$ to have a special 
rational point, either $\QQ (\sqrt{-p})$ or $\QQ(\sqrt{-251})$ must have class 
number one or the field $\QQ (\sqrt{-251p})$ must have class number two (see 
e.g. \cite{RotSkoYaf}, section 4). These conditions are not satisfied by the 
primes $p$ satisfying our conditions. 

In general it is difficult to construct infinite families of counter-examples to Hasse 
principle among the curves of this form (see \cite{RotSkoYaf}). However, using 
a PARI computer program, one can determine, for a given example, whether the 
curve has local points everywhere. We have done this for a particular example of 
$p=137$. This prime satisfies the congruences above. Note however that, due to 
the ineffectiveness of Lemma 3.5, we cannot guarantee the smoothness of the 
curve $\widetilde{X^{137 \cdot 251}/w_{251}}$ along the exceptional 
$\FF_{137}$-rational component. Using the theorem 3.1 of \cite{RotSkoYaf}, 
we verify first that the curve $X^{137\cdot 251}/w_{251}$ has rational points 
over $\RR$, $\QQ_{137}$ and $\QQ_{251}$. This is implied by the conditions 
\begin{enumerate}
\item $\left( \frac{251}{137}\right)=-1$.
\item $251 \equiv 3 \mod 4$ and $137 \equiv 1 \mod 4$.
\end{enumerate}
Then one verifies, that for every $\ell < 4g^2$ ($g$ being the genus
of the curve) one of the quantities (notations of \cite{RotSkoYaf}) $\Sigma_\ell 
(251\cdot 137)$ or $\Sigma_{251\ell} (251\cdot 137)$ is non-zero. This
guarantees the existence of local points at all primes of good reduction (i.e. $\ell 
\not=137,251$). It becomes difficult to verify this condition for primes $p> 137$ 
as one exceeds the PARI limit on precompiled primes.

We state our theorem for an arbitrary $p$.
\begin{thm}
For every large enough prime $p$ satisfying the congruences of Lemma 6.1 
the curve $X^{p\cdot 251}/w_{251}$ does not have $\QQ$-rational points.
\end{thm}

%
%
%

%
\begin{picture}(10,210)
\put(0,195){\line(20,0){90}}
\put(0,180){Pierre Parent}
\put(0,165){Institut de Math\'ematiques de Bordeaux}
\put(0,150){Universit\'e de Bordeaux I}
\put(0,135){351, cours de la lib\'eration}
\put(0,120){33 405 Talence C\'edex France}
\put(0,105){Adresse \'electronique : {\tt Pierre.Parent@math.u-bordeaux1.fr}}
\put(0,90){ }
\put(0,75){Andrei Yafaev}
\put(0,60){Department of Mathematics}
\put(0,45){University College London}
\put(0,30){Gower street}
\put(0,15){London, WC1E 6BT England}
\put(0,0){E-mail adress: {\tt yafaev@math.ucl.ac.uk}}

\end{picture}

\end{document}